\documentclass[
dvipdfm,
a4paper,reqno,12pt]{amsproc}%

\usepackage[pdftitle={Moments},
        verbose,
            hyperindex=true,
             naturalnames=true,
              linkcolor=red]{hyperref}
\usepackage{a4}

\usepackage{fullpage}
\usepackage{amsthm}
\usepackage{amsmath}

{
\theoremstyle{plain}
\newtheorem{theorem}{Theorem}
\newtheorem{lemma}{Lemma}
\newtheorem{corollary}{Corollary}
}

{ \theoremstyle{definition}
}

\theoremstyle{plain}

\theoremstyle{remark}

\def\Com#1{\mathbb{C}^{#1}}
\def\R#1{\mathbb{R}^{#1}}
\def\DD{\mathbb{D}}
\def\Pol#1{{\mathcal P}_{#1}(\overline{\DD})}
\def\Polc#1{{\mathcal P}_{#1,\,\, loc}(\overline{\DD})}

\def\Hm{\mathcal{T}}
\def\Bm{\mathcal{B}}

\def\z{\zeta}
\def\fr#1{\frac{1}{#1}}
\DeclareMathOperator*{\res}{res}
\DeclareMathOperator*{\Res}{\mathcal{R}}

\begin{document}

\title[Complex moment mapping]{Ullemar's formula for the Jacobian \\of the complex moment mapping}

\author{Kuznetsova O.S.}
\address{Volgograd State University, 2-ya Prodolnaya 30\\ Volgograd, 400062
Russia}

\curraddr{Royal Institute of Technology, Lindstedtsv\"agen 25,
10044 Stockholm, Sweden} \email{astra1987@mail.ru}

\author{Tkachev V.G.}
\address{Volgograd State University, 2-ya Prodolnaya 30\\ Volgograd, 400062
Russia}

\curraddr{Royal Institute of Technology, Lindstedtsv\"agen 25,
10044 Stockholm, Sweden} \email{Vladimir.Tkachev@volsu.ru}

\thanks{The paper is supported by the G\"oran Gustafsson
Foundation and the Russian Foundation for Basic Research, grant
03-01-00304}



\date{\today}

\begin{abstract}
The complex moment sequence $\mu(P)$ is assigned to a univalent
polynomial $P(z)$ by the Cauchy transform of the domain $P(\DD)$,
where $\DD$ is the unit disk. We establish the representation of
the Jacobian $\det d\mu(P)$  in terms of roots of the derivative
$P'(z)$. Combining this result with the special decomposition for
the Hurwitz determinants, we prove a formula for $\det d\mu(P)$,
which was previously conjectured by C.~Ullemar. As a consequence,
we show that the boundary of the class of all locally univalent
polynomials in $U$ is contained in the union of three irreducible
algebraic surfaces.
\end{abstract}

\maketitle
\section{Introduction}

Let $f(z)$, $f(0)=0$, be an analytic function defined in a
neighborhood of the unit disk $\DD$ and $k\geq 0$ be a nonnegative
integer. Then the \textit{complex moments} of $f(z)$  are defined
by
\begin{equation*}\label{eq1} M_k(f)=\frac{i}{2\pi
}\iint\limits_{\DD}f^k(z) |f'(z)|^2 \,dz\wedge d\overline{z}=
\frac{i}{2\pi }\iint\limits_{D}\z^k \,d\z\wedge d\overline{\z}.
\end{equation*}

This notion appears in several problems of complex analysis and
its applications. In particular, if $f(z)$ is a univalent function
in $\DD$, then the latter sequence constitutes an infinite family
of invariants of the Hele-Shaw problem \cite{R72}. On the other
hand, the sequence $(M_k(f))_{k\geq 0}$ defines the germ at
infinity of the Cauchy transform of the domain $D$
$$
\hat{\chi}_{D}=\frac{1}{\pi}\iint_{D}\frac{d\zeta}{z-\zeta}=\sum_{k\geq
0}\frac{M_k(f)}{z^{n+1}}.
$$

Since the above definition may be regarded as a two-dimensional
extension of the Stieltjes moments on the real line \cite{Ak}, it
makes natural the corresponding inverse problem of defining $f(z)$
by its moment sequence. It follows from the result of M.~Sakai
\cite{Sak} that without any additional restrictions, $f(z)$ (or
the domain $D$) can not be uniquely determined by its moments.
Some recent results concerning the reconstruction of a domain by
its complex moments can be found in \cite{GSS},
\cite{GP},\cite{MVCW}.

Throughout this paper we suppose that $f(z)$ is a polynomial
\begin{equation}\label{p}
P(z)=a_1z+\ldots+a_nz^n, \qquad a_1>0,
\end{equation}
of degree $n\geq 2$. Then $P(\DD)$ is an example of quadrature
domain (see \cite{Ah} and \cite[p.~11]{EtVarch1}). It follows from
a formula of S.~Richardson (see (\ref{summa}) below) that in this
case the corresponding sequence $(M_k(P))_{k\geq 0}$ is finite and
\begin{equation}\label{mzero}
M_k(P)=0,  \qquad k\geq n=\deg P.
\end{equation}
Moreover,
$$
M_0(P)=\sum_{j=1}^{n}j|a_j|^2>0, \qquad
M_{n-1}(P)=a_1^{n}\overline{a}_n\ne 0.
$$
Then  it follows from (\ref{mzero}) and Richardson's formula
(\ref{summa}) below that the complex moment sequence induces the
\textit{moment mapping} as a polynomial mapping
\begin{equation}\label{mu-com}
\mu_{\Com{}}( P) = (M_0(P),\ldots ,M_{n-1}(P)): \R{+}\times
\Com{n-1}\to\R{+}\times \Com{n-1}.
\end{equation}

Similarly, in the case of the \textit{real} polynomials $P(z)$,
i.e. $a_k\in\R{}$, (\ref{mu-com}) induces a real polynomial
mapping
\begin{equation}\label{mu}
\mu( P) = (M_0(P),\ldots ,M_{n-1}(P)): \R{+}\times
\R{n-1}\to\R{+}\times \R{n-1}.
\end{equation}

Thus the above mentioned inverse problem can be reformulated as an
injectivity problem for the preceding polynomial mappings.

Let $\Pol{n}$ denotes the class of all polynomials (\ref{p})
univalent in a neighborhood of the closed unit disk, $a_k\in\R{}$,
$1\leq k\leq n$. By $\Polc{n}$ we denote the class of the locally
univalent polynomials; equivalently, $P(z)\in \Polc{n}$ iff
$P'(z)\ne 0$ in $\overline{\DD }$. It is clear, that $\Pol{n}$ is
a proper subclass of $\Polc{n}$ for $n\geq 3$.

The main difficulty in the study of the injectivity problem of
$\mu$ and $\mu_{\Com{}}$ is a highly involved structure of the
class of univalent polynomials $\Pol{n}$. Only some low degree
($n\leq 3$) results are known (see \cite{Ks}, \cite{Kw},
\cite{Suf}).

It was proven by C.~Ullemar in \cite{Ul} that $\mu $ is
\textit{globally} injective on $\Pol{3}$ and the injectivity
property fails on $\Polc{3}$. The first general result for the
\textit{locally }univalent polynomials (actually, even with
complex coefficients) is due to B.~Gustafsson \cite{BG84} and
states that $\mu$ is locally injective on $\Polc{n}$. The question
whether $\mu$ is \textit{globally} injective on $\Pol{n}$ for
$n\geq 4$ is still open.

In her paper, Ullemar conjectured the following formula for the
Jacobian of $\mu$:
\begin{equation}\label{J}
J(P)\equiv\det \, d\mu
(P)=2{}^{-\frac{n(n-3)}{2}}a_1{}^{\frac{n(n-1)}{2}}P'(1)P'(-1)\Delta
_n(\widetilde{P'}(z)),
\end{equation}
which will be in focus of the present paper. Here
$\Delta_n(\widetilde{P'}(z))$ denotes the main Hurwitz determinant
for the M\"obius transformation of the derivative $P'(\zeta)$ (see
exact definitions in section~\ref{SecGurv}).


A useful feature of (\ref{J}) is that it immediately implies the
local injectivity property. Indeed, by the well-known Hurwitz
theorem  the inner determinant in (\ref{J}) is positive when
$P'(z)$ has no roots in a right half plane.

Our first result gives the following alternative formula for
evaluation of $J(P)$ via the inner characteristics of $P$.

\begin{theorem}[Derivative Roots Formula]
Let $P(z)=a_1z+\ldots +a_nz^n$, $a_k \in \R{}$ and $\zeta _1$,
$\ldots$, $\zeta _{n-1}$ are all zeroes of the derivative $P'(z)$.
Then
\begin{equation}
\begin{split}
 J(P)&= 2a_1^{\frac{n(n-1)}{2}} (na_n)^n\cdot
 \prod_{i\leq j}(\zeta _i\zeta _j-1)=\\
        &= 2a_1^{\frac{n(n-1)}{2}}(na_n)^{n-2} \,P'(1)\,P'(-1)\prod_{i<j}(\zeta
_i\zeta _j-1)
\end{split}
 \label{Jac}
\end{equation}
\label{th1}
\end{theorem}

Actually, the right hand side of (\ref{Jac}), as a symmetric
function of the roots, can be represented as a homogeneous form
$$
J(P)=2a_1^{\frac{n(n-1)}{2}} V_{n-1}(b_1,\ldots,b_n)\sum_{j=1}^n
b_j\sum_{k=1}^n (-1)^kb_k,
$$
where $b_k=ka_k$ are the coefficients of $P'(z)$ and $V_{n-1}$ is
a homogeneous irreducible polynomial of degree $(n-1)$ (see
section~\ref{sec:5} for precise definitions).

\begin{theorem}[Resultant Formula]
Let $A^*(z)=z^p A(1/z)$ be the reciprocal polynomial to
$A(z)=\alpha_0+\alpha_1z+\ldots +\alpha_pz^p$. Then
\begin{equation*}\label{-A31}
J(P)^2=4(-1)^{n-1}a_1^{n(n-1)}\Res (P',{P'}^*)\cdot P'(-1)P'(1),
\end{equation*}
where $\Res (A,B)$ denotes the resultant of the corresponding
polynomials. \label{theo:resul}
\end{theorem}

Now,  the Ullemar formula (\ref{J}) can be obtained as a
consequence of Theorem~\ref{th1} and certain auxiliary properties
of the Hurwitz determinants which we get in section~\ref{SecGurv}.

As another application we give an alternative proof of the above
mentioned result of Gustafsson.

\begin{corollary}
\label{corGus} The mapping $\mu(P)$ is locally injective on the
set $\Polc{n}$, $n\geq 1$.
\end{corollary}

\begin{proof}
Indeed, given any polynomial $P(z)\in \Polc{n}$ with real
coefficients we have $a_n\ne 0$ and $a_1=P'(0)\ne 0$. Moreover,
$|P'(\zeta)|\ne 0$ in $\overline {\DD }$ and it follows that all
zeroes of the first derivative $|\zeta _k|>1$, $k=1,\ldots ,n-1$.
Thus (\ref{Jac}) implies that $J(P)\ne 0$.
\end{proof}

It turns out that $J(P)$  can be associated with the structural
properties of the class $\Polc{n}$ as follows. Let us identify a
polynomial $P(z)=\sum_{j=1}^na_jz^j$ with  the point
$(a_1,\ldots,a_n)\in\R{n}$ and put
$$
\mathcal{P}^n_{loc}=\cup_{1\leq j\leq n}\Polc{j}.
$$

\begin{theorem}\label{ttt}
Let $n\geq 3$, then the boundary  of $\mathcal{P}^n_{loc}$ is
contained in the union of the following three irreducible
algebraic varieties: the hyperplanes
\begin{equation}\label{Pi}
\begin{split}
&\Pi^{+}:\quad P'(1)=a_1+ 2a_2+\ldots +na_n=0,\\
&\Pi^{-}:\quad  P'(-1)=a_1- 2a_2+\ldots +(-1)^{n-1}na_n=0,
\end{split}\end{equation}
and an algebraic surface of $(n-1)$th order given by
\begin{equation}\label{Pi1}
\mathcal{A}:\quad  V_{n-1}(a_1,2a_2,\ldots, na_n)=0.
\end{equation}
\end{theorem}

It follows from the preceding results that $\mathcal{P}^n_{loc}$
is exactly an open component of the set $\{P:J(P)\ne 0\}$.

The similar result for the univalent classes $\Pol{n}$ is due to
Quine \cite{Q}. But only upper estimates for the degree of the
boundary $\partial \Pol{n}$ have been established there.

We notice that the previous formulae as well as the suitable
modifications of basic facts below are still valid for polynomials
with  complex coefficients. This will be accomplished in a
forthcoming paper.

\medskip
\textbf{Acknowledges:} The authors wish to thank Bj\"orn
Gustafsson for bringing their attention to the present theme and
for fruitful discussions. We are grateful to Harold Shapiro for
his helpful comments and the referee for careful reading and
suggestions which led to an improvement of the paper.

\section{Preliminary results}
\label{subsec:21}

Following to  S.~Richardson \cite{R72} one can write the following
expressions for $M_k(P)$
\begin{equation}
M_k(P)=\sum i_1 a_{i_1}\cdots a_{i_{k+1}}\overline{ a}_{i_1+\ldots
+i_{k+1}}, \label{summa}
\end{equation}
where the sum is taken over all possible sets of indices $i_1$,
$\ldots $, $i_k\geq 1$. It is assumed that $a_j=0$ for $j\geq
n+1$. These formulae are easy to use for straightforward
manipulations with the complex moments and it follows that
$\mu_{\Com{}}$ as well as $\mu$ are polynomial mappings.
Nevertheless, this representation is  useless for the further
study of analytic properties of $\mu$.

We shall use in the sequel the following simple residue
representation of the moment sequence for real polynomials
\begin{equation}
M_k(P)=\frac{1}{k+1}\res_{\zeta=0}\left(P^{k+1}(\zeta)
P'\left(\frac{1}{\zeta}\right)\frac{1}{\zeta}\right). \label{res}
\end{equation}

Indeed, it follows from Stokes' formula that
\begin{equation}
\frac{i}{2\pi}\iint\limits_{G}w^k\, dw\wedge d\overline{
w}=\frac{i}{2\pi(k+1)}\int\limits_{\partial G}w^{k+1}\,
d\overline{w}, \label{Stm}
\end{equation}
where $G$ is an arbitrary $2$-chain in the complex plane. Letting
$G=P(\DD )$ and taking into account that $\bar \zeta =\zeta^{-1}$
on $\partial \DD $,  and the fact that
$\overline{P'(z)}=P'(\overline{z})$ for polynomials with real
coefficients, we obtain from (\ref{Stm})
$$
M_k(P)=\frac{i}{2\pi(k+1)}\int\limits_{\partial \DD
}P^{k+1}(\zeta)\overline{P'(\zeta)}\, d\bar \zeta= \frac{1}{2\pi
(k+1)}\int\limits_{\partial \DD
}P^{k+1}(\zeta)P'\left(\frac{1}{\zeta}\right)\,
\frac{d\zeta}{\zeta ^2}.
$$
This proves (\ref{res}).

Moreover, since $P(0)=0$, it follows that $P(\zeta)=zP_1(z)$,
where $P_1$ is a polynomial. Thus, the expression
$$
\zeta ^{k+1}P_1^{k+1}(\zeta)P'\left(\frac{1}{\zeta}\right)=\zeta
^{k-n}(a_1+\ldots +a_n\zeta ^{n-1})(a_1\zeta ^{n-1}+\ldots +a_n)
$$
is also a polynomial for all $k\geq n$ and it follows from
(\ref{res}) that
$$
\label{resM} M_k(P)=\frac{1}{k+1}\res_{\zeta=0}\zeta ^{k+1}
P_1^{k+1}(\zeta) P'\left(\frac{1}{\zeta}\right)=0
$$
which proves (\ref{mzero}). Therefore, the mapping $\mu$ in
(\ref{mu}) is well-defined.

Given two meromorphic functions $H_1$ and $H_2$ we write
$$
H_1(z)\equiv H_2(z) \mod [m_1;m_2]
$$
if the Laurent series of $H_2-H_1$   does not contain $z^m$ with
$m_1\leq m\leq m_2$.
\begin{lemma}
For any $k$, $0\leq k\leq n-1$,
\begin{equation}
P'(z)\left(P^k(z)+P^k\left( \frac{1}{z}\right)\right)\equiv
\sum_{\nu =1}^{n}\frac{\partial M_k(P)}{\partial a_{\nu}}\cdot
z^{\nu -1} \mod  [0;n-1]. \label{mod}
\end{equation}
\label{l1}
\end{lemma}

\begin{proof}
Let $ \lambda _m(f(z))=\res_{z=0}(f(z)z^{-1-m})$; then it follows
from the relations
$$
\frac{\partial P(1/z)}{\partial a_{\nu}}=\frac{1}{z^{\nu }},
\qquad \frac{\partial P'(z)}{\partial a_ {\nu}}=\nu z^{\nu -1},
$$
and ~(\ref{res}) that
\begin{equation}
\frac{\partial M_k(P)}{\partial a_ {\nu}}=\lambda
_0\left(P^{k}(1/z)P'(z)z^{1-\nu }\right)+\frac{\nu}{k+1}\lambda_0
\left(P^{k+1}(1/z){z}^{\nu }\right). \label{partM}
\end{equation}
On the other hand, integrating by parts yields
$$
\lambda _0\left( z^{\nu}P^{k+1}(1/z)\right)=\frac{1}{2\pi i}\int
\limits _{\partial \DD }P^{k+1}\left(1/z\right)z^{\nu
-1}\,dz=\frac{1}{2\pi i\nu }\int \limits _{\partial \DD
}d\left(z^{\nu }P^{k+1}\left(1/z\right)\right)+
$$
$$
+\frac{k+1}{2\pi i\nu}\int \limits _{\partial \DD
}P^{k+1}\left(1/z\right)P'\left(1/z\right) z^{\nu -2}\,dz=
\frac{k+1}{\nu}\lambda
_0\left(P^{k}\left(1/z\right)P'\left(1/z\right) z^{\nu-1 }\right),
$$
and taking into account that $\lambda_0(f(1/z))=\lambda _0(f(z))$
we arrive at
\begin{equation}\label{ppp}
    \lambda _0\left(
P^{k+1}\left(z^{-1}\right)z^{\nu}\right)=\frac{k+1}{\nu}\lambda
_0\left({P^k(z)P'(z)}{z^{1-\nu}}\right).
\end{equation}
Combining (\ref{ppp}) and (\ref{partM}), we get
\begin{equation*}\begin{split}
 \frac{\partial M_k(P)}{\partial a_ {\nu}} &=
\lambda
_0\left[P'(z)z^{1-\nu}\biggl(P^k(z)+P^k\left(z^{-1}\right)\biggr)\right]
=\\
&=\lambda_{\nu-1}\left[P'(z)\biggl(P^k(z)+P^k\left(z^{-1}\right)\biggr)\right]
\end{split}
\end{equation*}
and the required formula (\ref{mod}) follows.

\end{proof}

We notice that for any index $k\in \{0, \ldots ,n-1\}$ the
following expansion
\begin{equation}
P^k(z)+P^k\left(z^{-1}\right)=\sum_{m=-nk}^{nk}h_m^{(k)}z^m,
\label{13}
\end{equation}
yields the symmetry property: $ h_m^{(k)}=h_{-m}^{(k)}$.

To study (\ref{mod}) it is convenient to consider a slightly more
general case. Namely, given an arbitrary vector
$x=(x_0,x_1,\ldots,x_{n-1})$,  we define the following Toeplitz
matrix
$$
\Hm(x)= \left(
\begin{array}{cccc}
x_0 & x_1 & \cdots & x_{n-1}\\
x_1 & x_0 & \cdots & x_{n-2}\\
\vdots & \vdots & \ddots &\vdots \\
x_{n-1} & \cdots & \cdots & x_0
\end{array} \right).
$$
Then we can introduce the \textit{dual} matrix $\Bm (y)$,
$y=(y_0,y_1,\ldots,y_{n-1})$,  by
\begin{equation}
\Hm (x)\cdot y^\top= \Bm (y)\cdot x^\top, \qquad \forall
x\in\R{n}. \label{H}
\end{equation}
Unlike $\Hm (x)$, the matrix $ \Bm (y)$ is not symmetric and has a
more complicated structure. We shall study  $B(y)$ in more detail
the next section.

Let now $H_k(z)$ be rational functions having Laurent series of
the form
$$
H_k(z)=\sum_{m=-N}^{N}h_{|m|}^{(k)}\;z^m,
$$
and let
$$
B(z)=b_0+b_1z+\ldots +b_{n-1}z^{n-1}
$$
be a polynomial such that $b_{n-1}\ne 0$.

Then we can define polynomials
$$
\Phi _k(z)=\sum_{\nu =0}^{n-1}\varphi ^{(k)}_{\nu }z^{\nu}, \qquad
0\leq k\leq n-1
$$
such that
\begin{equation}
B(z)\cdot H_k(z)\equiv \Phi _k(z) \mod  [0;n-1]. \label{Qmod}
\end{equation}

Further, we consider the vectors $h^{(k)}=(h_0^{(k)},\ldots
,h_{n-1}^{(k)})$ and $b=(b_0,\ldots ,b_{n-1})$. It follows then
from (\ref{Qmod}) that the following matrix identity holds
$$
(\varphi_0^{(k)},\ldots ,\varphi_{n-1}^{(k)})^\top
\equiv{\varphi^{(k)}}^\top =\Hm (h^{(k)})\cdot b^\top,
$$
which by virtue of (\ref{H}) implies $ {\varphi ^{(k)}}^\top= \Bm
(b)\cdot {h^{(k)}}^\top$, $0\leq k\leq n-1. $ Therefore, denoting
by $\Phi$ and $H$ the matrices formed by combination of the
columns ${\varphi ^{(k)}}^\top$  and ${h^{(k)}}^\top$
respectively, we get $\Phi= \Bm(b) H$ and
\begin{equation}
\det \,\Phi=\det \, \Bm(b)\cdot \det \, H. \label{detbc}
\end{equation}

To apply the preceding arguments to our case we let $B(z)=P'(z)$
and $H_k(z)=P^k(z)+P^k(1/z)$. Hence, we obtain from (\ref{mod})
\begin{equation}\label{e:phi}
\varphi_{\nu}^{(k)}=\frac{\partial
M_k(P)}{\partial a_\nu}, \qquad d\mu (P)=\Phi.
\end{equation}

Thus, the problem of evaluating the Jacobian $J(P)$ can be
reduced, by virtue of (\ref{detbc}), to the corresponding problem
for the determinants of $\Bm(b)$ and $H$ (here $b_{j-1}=ja_j$
corresponds to the coefficients of $P'(z)$).

The latter determinant can be found as follows. First note that
$\|h_i^{(k)}\|$ is a lower triangular matrix in our case. Indeed,
we have $P(z)=zP_1(z)$, where $P_1(z)$ is a polynomial, and it
follows that
$$
P^k(z)+P^k\left(z^{-1}\right)=z^kP_1(z)+\frac{1}{z^k}P_1^k\left(z^{-1}\right)
=\sum_{m=k}^{kn}(z^m+z^{-m})h_m^{(k)}.
$$
This representation easily implies that $h_m^{(k)}=0$, where
$0\leq m\leq k-1$. Moreover, we have for the diagonal elements
$h_0^{(0)}=2$ and $h_k^{(k)}=a_1^k$. This yields
\begin{equation}
\det \, H=\det \,\|h_i^{(k)}\|=2\cdot a_1\cdot a_1^2\cdot \ldots
\cdot a_1^{n-1}= 2a_1^{\frac{n(n-1)}{2}}. \label{findet}
\end{equation}

\section{Toeplitz determinants}\label{PROOF}
The explicit expression of $\det \, \Bm (y)$ in terms of
the coefficients $y_0,\ldots ,y_{m}$ is messy and useless for the
further analysis. However, it turns out that this determinant can
be easily written in terms of certain intrinsic characteristics of
$y$. Namely, let us associate with any vector $y\in\R{m}$ the
polynomial
$$
B_y(z)=y_0+y_1z+\ldots +y_{m}z^{m},\qquad y_{m}\ne 0.
$$

\begin{theorem}
Let $\zeta_1$, $\ldots$, $\zeta_m$ be the roots of $B_y(\zeta)$
counted according to their multiplicities. Then
\begin{equation}
\det \, \Bm (y)=y^{m+1}_{m}\prod _{i\geq j}(\zeta _i\zeta _j-1),
\label{detB}
\end{equation}
\label{th2}
\end{theorem}

\begin{proof}
First note that the left-hand side of ~(\ref{detB}) is an
algebraic function of $y_0$, $\ldots$, $y_{m}$ and, hence, it is
sufficient to prove (\ref{detB}) for any $\zeta =(\zeta _1,\ldots
,\zeta _{m})$ outside a proper algebraic submanifold of $\Com{m}$.
Namely, we will suppose that $\zeta _i\ne \zeta _j$ for $i\ne j$
and $\zeta _i\zeta _j\ne 1$ for all $i$, $j$.

Given a nonnegative integer $k$ and $\zeta \in \Com{}$ we define
the following vector
$$
\{\zeta \}_{k}=(0,\ldots,0,1,\zeta ,\zeta ^2,\ldots ,\zeta
^{m-k})^\top\in\Com{m+1}, \qquad \{\zeta \}\equiv \{\zeta \}_{0}.
$$
Then letting $x=\{\zeta \}^\top$ in (\ref{H}) we get
$$
\Hm (\{\zeta \})\cdot y^\top = B_y(\zeta) \cdot \{\zeta^{-1} \}+
\sum_{i=0}^{m-1}y_i\left(\{\zeta \}_i-\{\zeta^{-1} \}_i\right)
$$
and changing $\zeta$ by $\zeta ^{-1}$ in the preceding formula we
arrive after summation at
\begin{equation}
\Hm (\{\zeta \}+\{\zeta ^{-1}\})\cdot y^\top=B_y(\zeta) \cdot
\{\zeta ^{-1}\} +B_y(\zeta^{-1}) \cdot \{\zeta\}. \label{HQb}
\end{equation}

Let $\zeta =\zeta _i$ be a root of $B_y(\zeta)$; then it follows
from~(\ref{HQb}) that
\begin{equation}
\Hm (\{\zeta _i\}+\{\zeta_i^{-1}\})\cdot y^\top=B_y(\zeta_i^{-1})
\cdot \{\zeta _i\} \label{zet}
\end{equation}
and
\begin{equation}
\Hm (e)\cdot y^\top=B_y(1)\cdot e^\top, \label{HQe}
\end{equation}
where $e=(2,\ldots ,2)\in\Com{m+1}$. Applying (\ref{H}) to the
left-hand sides of (\ref{zet}) and (\ref{HQe}) we obtain
\begin{equation*}
 \Bm (y)(\{\zeta _i\}+\{\zeta_i^{-1}\})=B_y(\zeta_i^{-1}) \cdot
\{\zeta _i\}, \qquad \forall  i=1,\ldots ,m,\label{B}
\end{equation*}
and
\begin{equation*}
 \Bm (y)\cdot e^\top=B_y(1)\cdot e^\top, \label{Be}
\end{equation*}

Combining the preceding expressions  into the matrix form we
arrive at the following relation for determinants
\begin{equation}
\begin{split}
   \det\Bm (y) &\det \mathcal{W}(1,\zeta _1, \ldots
,\zeta_{m})=\\
&=2B_y(1)
  \det \begin{pmatrix}
    1 & 1 &\ldots & 1\\
    1 & \zeta_1 &\ldots & \zeta_1^m\\
    \vdots& \vdots&\ddots&\vdots\\
    1 & \zeta_m &\ldots & \zeta_m^m
  \end{pmatrix}\prod_{j=1}^{m}B_y(\zeta_j^{-1})=\\
  &=2 (-1)^mB_y(1)
\prod_{k=1}^{m}B_y(\zeta_k^{-1})\cdot \prod_{i<j}(\zeta
_j-\zeta_i)\cdot \prod _{i=1}^{m}(1-\zeta _i),\label{BW}
\end{split}
\end{equation}
where $\mathcal{W}(\alpha _0,\alpha_1,\ldots ,\alpha _m)$ denotes
the matrix with the entries
 $$
 \mathcal{W}_{ij}=\|\alpha
_j^{i}+\alpha _j^{-i}\|_{i,j=0}^{m}.
$$

The determinant of $\mathcal{W}(\alpha _0, \ldots ,\alpha _m)$ can
be found by the same method as the Vandermonde determinant (see
also \cite[Part~ 4]{Det}):
$$
\det \,\mathcal{W}(\alpha _0, \ldots ,\alpha
_m)=\frac{2}{(\alpha_0 \ldots \alpha _m)^{m}}\prod_{i<j}(\alpha
_j-\alpha_i)\prod _{i<j}(\alpha_i \alpha _j-1)
$$
and it follows that
\begin{equation}\label{detW}
\det \,\mathcal{W}(1, \zeta_1,\ldots ,\zeta _m)=
\frac{2}{(\z_1\cdots \z_m)^m}\prod_{i<j}(\z _j-\z_i) \prod
_{i<j}(\z_i \z _j-1)\prod_{j=1}^m (1-\z_j)^2.
\end{equation}
On the other hand,
$$
\prod_{k=1}^{m}B_y(\zeta _k^{-1})=(-1)^m\frac{y_{m}^m}{(\z_1\cdots
\z_m)^m}\prod_{j=1}^{m}\prod_{i=1}^{m}(\zeta _i\zeta _j-1).
$$

Thus, applying the previous identities to (\ref{BW}) we obtain
$$
 \det \, \Bm (y)=\frac{B_y(1)y_m^m }{\prod_{i=1}^{m}(1-\zeta _i)}\prod_{j\geq i}(\zeta _i\zeta _j-1)
$$
which implies by virtue of
$$
B_y(1)=y_m\prod_{i=1}^{m}(1-\zeta _i)
$$
the required identity.
\end{proof}

\begin{proof}[Proof of Theorem~\ref{th1}] It follows from (\ref{e:phi}) that
 $d\mu (P)=\Phi$. Then applying (\ref{findet}) and Theorem ~\ref{th2} to
(\ref{detbc}) we obtain
$$
J(P)\equiv \det \, \left[\frac{\partial M_k(P)}{\partial a_i}
\right]=2a_1^{\frac{n(n-1)}{2}}\cdot b_{n-1}^n \prod_{i\leq
j}(\zeta _i\zeta _j-1),
$$
where $b_{n-1}=na_n$ is the leading coefficient of $B(z)\equiv
P'(z)$ and the theorem follows.
\end{proof}

\section{Hurwitz determinants and Ullemar's formula}\label{SecGurv}

Let us consider an arbitrary polynomial $R(z)=r_0+r_1z+\ldots
+r_mz^m$ of  degree $m\geq 1$. Let us extend the sequence of the
coefficients $r_k$ such that $r_k=0$ for all $k>m$ and $k<0$. Then
the  $m\times m$-matrix
\begin{equation*}\label{meb2}
\mathcal{G}(R)\equiv
\begin{pmatrix}
    r_{m-1} & r_{m-3}&\ldots& r_{1-m} \\
    r_{m} & r_{m-2}& \ldots&r_{2-m} \\
    \vdots& \vdots&\ddots&\vdots\\
    r_{2m-2}&r_{2m-4}&\ldots&r_{0}
  \end{pmatrix}
\end{equation*}
is called the Hurwitz matrix  of the polynomial $R(z)$
\cite{Gant}. More specifically, the entries of the matrix have the
form
\begin{equation}\label{meb3}
\mathcal{G}_{ij}(R)=r_{m+i-2j}.
\end{equation}

The main diagonal minor $\Delta(R)$ of the $(m-1)$th order of
$\mathcal{G}(R)$ is said to be the \textit{Hurwitz determinant} of
$R$. It immediately follows from the above definition that
\begin{equation}\label{minor}
  \det \mathcal{G}(R)=r_0 \;\Delta(R).
\end{equation}

\begin{theorem}\label{th3}
The Hurwitz determinant of $R(z)$, $\deg R=m$, has the following
representation
\begin{equation}\label{meb5}
\Delta(R)=(-1)^{\frac{m^2-m}{2}}r_m^{m-1}\prod_{1\leq i < j \leq
m}(z_i+z_j),
\end{equation}
where $z_i$ are the roots of $R(z)$ counted according to their
multiplicity.
\end{theorem}

Before we give the proof of the theorem let us formulate some of
its corollaries. Let us consider the \textit{M\"obius
transformation}  of the polynomial $R(z)$ given by
\begin{equation*}\label{meb1}
\widetilde{R}(z)=(z+1)^m R\left(\frac{z-1}{z+1}\right)\equiv
\widetilde{r}_0+\widetilde{r}_1z+\ldots +\widetilde{r}_mz^m.
\end{equation*}
Obviously, $\zeta_k=\frac{1+z_k}{1-z_k}$ are the roots of
$\widetilde{R}(z)$ whenever $z_1$, \ldots $z_m$ are the roots of
$R(z)$. In particular, all the roots of $R(z)$ are contained in
the unit disk if and only if the roots of $\widetilde{R}(z)$ lie
in the right half-plane. Moreover,
$$
\prod_{1\leq i < j \leq m}(\zeta_i+\zeta_j)=2^{\frac{m(m-1)}{2}}
\prod_{1\leq i < j \leq m}(1-z_iz_j)
\left(\prod_{i=1}^{m}(1-z_i)\right)^{1-m}.
$$

\medskip
\noindent Then the following identities
$$
\prod_{i=1}^{m}(1-z_i)=\frac{R(1)}{r_m},
$$
$$
\widetilde{r}_m=\lim_{z\to \infty} z^{-m}\widetilde{R}(z)=R(1),
$$
together with (\ref{meb5}) yield

\begin{corollary}\label{meb6}
In the previous notations
\begin{equation}\label{meb7}
\Delta(\widetilde{R})= 2^{\frac{m^2-m}{2}}r_m^{m-1}\prod_{1\leq i
< j \leq m}(z_iz_j-1),
\end{equation}
where $r_m$ is the leading coefficient of $R$ and $\{z_i\}_{1\leq
i\leq m}$ are  the roots of $R$.
\end{corollary}

Now, Ullemar's conjectured  formula (\ref{J}) is a simple
consequence of (\ref{meb7}) and Theorem~\ref{th1}.

\begin{corollary}[Ullemar Formula]\label{Ullem}
The Jacobian of the complex moment mapping $\mu$ has the following
representation
$$
J(P)\equiv\det \, d\mu
(P)=2{}^{-\frac{n(n-3)}{2}}a_1{}^{\frac{n(n-1)}{2}}P'(1)P'(-1)\Delta
_n(\widetilde{P'}(z)),
$$
where $n=\deg P$.
\end{corollary}

\begin{proof}[Proof of Theorem \ref{th3}]
Similarly to the proof of Theorem~\ref{th1} we can  assume that
$R(z)$ has no multiple roots. Then we have from (\ref{meb3}) for
any $\zeta\in\Com{}$ and any index $i$, $1\leq i\leq m$, that
\begin{equation*}\label{g1}
\sum_{j=1}^{m}\mathcal{G}_{ij}(R)\zeta^{2m-2j}=
\sum_{j=1}^{m}r_{m+i-2j}\zeta^{2m-2j}=
\zeta^{m-i}\sum\nolimits^\star_k \zeta^k r_k,
\end{equation*}
where the $k$th index in the last sum has the form
$k(i,j)=m+i-2j$, $j=1,\ldots,m$.

Now suppose that $i\in[1,m]$ and let $\bar i=(m-i)$. It is clear
that $k$ takes  only the even (or only odd) values which range
between $(-\bar i)$ and $(2m-2-\bar i)$ with changing $j$ between
$1$ and $m$. Moreover, both $k$ and $\bar i$ have the same parity
and
$$
-\bar i\leq k(i,j)\leq 2m-2-\bar i, \qquad \bar i=0,1,\ldots, m-1.
$$

Hence, for each $i$ from $[1,m]$ the indices $k(i,j)$ take all the
values of $\bar i$ from interval $\{0,1,\ldots, n\}$ when
$j\in[0,m]$. By virtue of this property we conclude that
\begin{equation}\label{g2}
\sum_{j=1}^{m}\mathcal{G}_{ij}(R)\zeta^{2m-2j}=\zeta^{m-i}R^{[\,\bar
i\,]}(\zeta),
\end{equation}
where  $R^{[p]}(\zeta)$ is the even (or odd) part of $R(\zeta)$:
$$
R^{[p]}(\zeta)= \frac{1}{2}\biggl(R(\zeta)+(-1)^pR(-\zeta)\biggr).
$$

Let now $\zeta=\zeta_k$ be the $k$th root of $R(z)$. Taking into
account that
$$
R^{[p]}(\zeta_k)=
R(\zeta_k)-R^{[p+1]}(\zeta_k)=-R^{[p+1]}(\zeta_k),
$$
we obtain
$$
R^{[p]}(\zeta_k)=(-1)^pR^{\mathrm{ev}}(\zeta_k).
$$
Here $R^{\mathrm{ev}}=R^{[0]}$ is the even part of $R(z)$ and we
see from (\ref{g2}) that
$$
\sum_{j=1}^{m}\mathcal{G}_{ij}(R)\zeta_k^{2m-2j}=(-\zeta_k)^{m-i}R^{\mathrm{ev}}(\zeta_k).
$$

Combining the last identities for $k=1,2,\ldots,m$ into the
matrices we obtain for their determinants
\begin{equation}\label{g4}
  \det \mathcal{G}(R)\det
\mathcal{V}(\zeta_1^2,\ldots,\zeta_m^2)
  = (-1)^{\frac{m(m-1)}{2}}\det
\mathcal{V}(\zeta_1,\ldots,\zeta_m)
\prod_{k=1}^{m}R^{ev}(\zeta_k),
\end{equation}
where $\mathcal{V}(a_1,\ldots,a_m)=\|a_j^{k-1}\|_{j,k=1}^{m}$ is
the Vandermonde matrix.

On the other hand, we have for the even part
$$
R^{ev}(\zeta_k)=\frac{1}{2}R(\zeta_k)=
\frac{r_m}{2}\prod_{i=1}^{m}(\zeta_i+\zeta_k),
$$
and it follows from (\ref{g4}) that
$$
\det \mathcal{G}(R)\cdot\prod_{1\leq i < j \leq
m}(\zeta_j^2-\zeta^2_i)= \frac{(-1)^mr_m^m}{2^m}
\prod_{i,j=1}^{m}(\zeta_i+\zeta_j)\prod_{1\leq i < j \leq
m}(\zeta_i-\zeta_j).
$$

Hence, applying (\ref{minor}) we find
\begin{equation*}\label{prod}
 \Delta(R)= \frac{(-1)^{\frac{m^2+m}{2}}r_m^m}{2^m r_0}
\prod_{{1\leq i \leq j \leq m}}(\zeta_i+\zeta_j)
\end{equation*}
and rewriting the last product as
$$ \prod_{{1\leq i \leq j
\leq m}}(\zeta_i+\zeta_j)= \prod_{{i=1}}(2\zeta_i) \prod_{{1\leq i
< j \leq m}}(\zeta_i+\zeta_j)=\frac{(-2)^m r_0}{r_m}\prod_{{1\leq
i < j \leq m}}(\zeta_i+\zeta_j),
$$
we arrive at the required identity and the theorem is proved.
\end{proof}

\section{Representations via the resultants}

Recall that given two polynomials
$$
A(z)=A_n(z-\alpha_1)\ldots(z-\alpha_n)=A_0+A_1z+\cdots+A_nz^n
$$ and
$B(z)=B_n(z-\beta_1)\ldots(z-\beta_n) $ the  product
\begin{equation*}\label{rez}
    \Res(A,B)=A_n^nB_n^n\prod_{i,j=1}^{n}(\alpha_i-\beta_j)
\end{equation*}
is called  the \textit{resultant} of $A$ and $B$.

If $A(z)$ and $B(z)$ are the mutually reciprocal polynomials
$$
B(z)=z^n A(1/z)\equiv A^*(z),
$$
then $ B_{n-j}=A_j$, $j=0,\ldots n$ and  we have for their roots:
$ \beta_j=\frac{1}{\alpha_j}$. Then the corresponding resultant
can be rewritten in the matrix form
\begin{equation}\label{rez_defA}
\Res(A,A^*)=\det
\begin{pmatrix}
  A_0 & A_1 & \ldots & \ldots     &A_n     &              \\
      & A_0 & A_1    &\ldots&  \ldots      &A_n           \\
   &  &\ldots&\ldots &\ldots &\ldots&\ldots \\
  &   &     & A_{0}  & A_{1} & \ldots  &\ldots &A_n    \\
  A_n & A_{n-1} & \ldots & \ldots     &A_0     &               \\
      & A_n & A_{n-1}    &\ldots& \ldots       &A_0             \\
     &  &\ldots&\ldots &\ldots &\ldots&\ldots \\
  &   &     & A_{n}  & A_{n-1} & \ldots  &\ldots &A_0    \\
\end{pmatrix}.
\end{equation}

It is easy to see that the last matrix is of $2n$th order and with
$A_0$ as its diagonal elements. On the other hand,
\begin{equation*}\label{AA}
\Res(A,A^*)=A_n^nA_0^n\prod_{i,j=1}^{n}\left(\alpha_i-\fr{\alpha_j}\right)=
\frac{A_n^nA_0^n}{(\alpha_1\ldots\alpha_n)^{n}}
\prod_{i>j}^{n}\left(\alpha_i\alpha_j-1\right)^2
\prod_{i=1}^{n}(\alpha_i^2-1)
\end{equation*}
and by Vi\`ete's theorem
\begin{equation*}
\alpha_1\cdots\alpha_n=(-1)^{n}\frac{A_0}{A_n},\qquad A(1)A(-1)=
A_n^2\prod_{i=1}^n(\alpha_i^2-1),\label{A0}
\end{equation*}
 we conclude that
\begin{equation}\label{A3}
\begin{split}
\Res(A,A^*)&=(-1)^nA(-1)A(1)A_n^{2n-2}
\prod_{i>j}^{n}\left(\alpha_i\alpha_j-1\right)^2 =\\
&=\frac{(-1)^nA_n^{2n+2}}{A(1)A(-1)} \left[\prod_{i\geq
j}^{n}\left(\alpha_i\alpha_j-1\right)\right]^2.
\end{split}
\end{equation}

Thus, we have from (\ref{A3})
\begin{equation}\label{A31}
W_n^2(A)=(-1)^n\Res(A,A^*)A(-1)A(1)
\end{equation}
where
\begin{equation}
W_n(A)=A_n^{n+1} \prod_{i\leq j}(\alpha_i\alpha_j-1). \label{JacW}
\end{equation}

As an immediate consequence of its definition, $W_n(A)\equiv
W_n(A_0,A_1,\ldots, A_n)$ is a homogeneous form of order $n=\deg
A$. Moreover, it admits the following factorization
\begin{equation}
W_n(A)=A(-1)A(1)V_n(A), \qquad V_n(A)=A_n^{n-1} \prod_{i<
j}(\alpha_i\alpha_j-1),\label{A7}
\end{equation}
where $V_n(A)$ is a homogeneous form of degree $(\deg A-2)$.

On the other hand, it follows from
\begin{equation}\label{equ:def1}
V_n(A)=A_0^{n-1}
\prod_{i<j}\left(1-\frac{1}{\alpha_i\alpha_j}\right),
\end{equation}
that we have the recursion formula
\begin{equation*}\label{reduc}
V_n(A_0,A_1,\ldots, A_k,0,\ldots,0)=A_0^{n-k}
V_{k}(A_0,A_1,\ldots, A_k).
\end{equation*}

Here are the explicit expressions for $V_k$:
\begin{equation*}
\begin{split}
V_3(A)&=A_0^2-A_0A_2+A_1A_3-A_3^2\\
V_4(A)&= A_4(-A_1^2+A_3A_1+A_4^2-A_4A_2-A_0A_4+2A_0A_2-
A_0^2)+\\
        & +A_0(A_0^2-A_0A_2+A_1A_3-A_3^2).
\end{split}
\end{equation*}

\begin{proof}[Proof of Theorem~\ref{theo:resul}]
Substituting the derivative
$$
P'(z)=a_1+2a_2z+\ldots+na_nz^{n-1}\equiv
b_1+b_2z+\ldots+b_nz^{n-1},
$$
as $A(z)$ to (\ref{A3}) and (\ref{A31}) we obtain
\begin{equation}\label{A5}
\left[b_n^n\prod_{i\geq j}^{n-1}\left(\z_i\z_j-1\right)\right]^2=
(-1)^{n-1} \Res(P',P'^*)P'(-1)P'(1).
\end{equation}

Then comparing the last relations with the definition (\ref{JacW})
we arrive at the following formula
\begin{equation*}\label{A6}
W_{n-1}(P')^2= (-1)^{n-1} \Res(P',P'^*)P'(-1)P'(1).
\end{equation*}
Finally, combining the preceding identity with  (\ref{Jac}) we
attain the required representation of $J(P)$
\begin{equation*}
J^2(P)=4b_1^{n^2-n}W_{n-1}^2(P')=(-1)^{n-1}
\Res(P',P'^*)P'(-1)P'(1) \label{new}
\end{equation*}
which completes the proof.
\end{proof}

The following property of $V_k$ will be used in the next section.

\begin{theorem}\label{teo1}
 $V_n(A)\equiv V_n(A_0,A_1,\ldots, A_n)\in\Com{}[[A_0,A_1,\ldots,
A_n]]$ is an irreducible polynomial.
\end{theorem}

\begin{proof}
A simple analysis of the denominator of the right-hand side of
(\ref{equ:def1}) shows that $A_n$ can not be a divisor of
$V_n(A)$. On the other hand, we notice that $V_n(A)$ can be
represented as a symmetric polynomial function of the roots
$(\alpha_k)_{1\leq k\leq n}$ of $A(z)=0$.

Let $H_1(A)$ and $H_2(A)$ be two nontrivial (i.e., different from
the identical constants) divisors of $V_n(A)$. It is a consequence
of the homogeneity of $V_n(A)$ that both of $H_k(A)$ are
homogeneous too. Moreover, in our assumptions $h_k=\deg H_k\geq
1$.

By Vi\'ete's theorem
$$
A_k=A_n\sigma_k (\alpha_1,\ldots, \alpha_n)
$$
where $\sigma_k$ is $k$th symmetric function of $(\alpha_1,\ldots,
\alpha_n)$. Then substituting  the last expressions for $H_k(A)$
yields by virtue of the homogeneity of $H_k$ that
$$
H_k(A)=A_n^{h_k}Y_k(\alpha_1,\ldots, \alpha_n)
$$
where the $Y_k$, $k=1,2$, are polynomials in $\alpha_j$. On the
other hand, it follows from (\ref{A7}) that $h_1+h_2=n-1$ and the
each $Y_k$ must be a divisor of
$$
\prod_{i< j}(\alpha_i\alpha_j-1).
$$
But the last product consists of irreducible factors
$(\alpha_i\alpha_j-1)$ only. Moreover, if one
$(\alpha_i\alpha_j-1)$ occures in $Y_1$ as a divisor then by
symmetry the others have to be the divisors as well.

It follows that one of $Y_k$ contains none $\alpha_i$, i.e. it has
the form $A_n^p$. Thus, applying the remark in the beginning of
the proof we see that $p=0$. But this means that $Y_k$ must be a
constant factor that contradicts to our assumption and proves the
theorem.
\end{proof}

\section{Proof of Theorem~\ref{ttt}}
\label{sec:5}

Let $P(z)=a_1z+\ldots a_nz^n$, $P\in \mathcal{P}^n_{loc}$, be a
locally injective in the unit disk polynomial. We identify $P'(z)$
with the vector of its coefficients
$$
b=(b_1,b_2,\ldots,b_n)\in\R{n},
$$
where $b_k=ka_k$. We write also $\Res (p,q)=\Res(P',Q')$ for the
corresponding vectors $p$ and $q$. Moreover,  by $S$ we denote the
differential operator regarded as a linear transform in $\R{n}$:
$$
S(P)=P'(z):\R{n}\to\R{n}.
$$

Then the following consequence of (\ref{A31}) and (\ref{A7})
\begin{equation}\label{sign}
b\in \ker W_{n-1} \Leftrightarrow \Res(b,b^*)=0
\end{equation}
is useful. Here  $b^*=(b_n,\ldots,b_1)$ corresponds to $P'^*$.

\begin{lemma}\label{lem:com}
The set $\mathcal{P}^n_{loc}$ is an open connected subset of
$\R{n}$. A polynomial $P(z)$ is an element of the boundary
$\partial \mathcal{P}^n_{loc}$ if and only if the following
properties hold

(i) $P'(z)$ contains no zeroes in $\DD$;

(ii) $\Res (P', P'^*)=0$.
\end{lemma}

\begin{proof}
The openness of $\mathcal{P}^n_{loc}$  obviously follows from the
fact that
\begin{equation}\label{www}
\min_{z\in\overline{\DD}}|P'(z)|>0,\qquad \forall P\in
\mathcal{P}^n_{loc}.
\end{equation}
Furthermore, let $P\in \mathcal{P}^n_{loc}$; then the homotopy
$$
a_\lambda=(a_1,a_2t,\ldots, a_nt^{n-1}), \qquad t\in[0;1],
$$
corresponds to the dilatation  $P_t(z)= \frac{1}{t}P(t z)$ and
connects $P(z)$ and $Q(z)=a_1z$
 inside of $\mathcal{P}^n_{loc}$ since $P_t'(z)=P'(t
z)\ne0$ in $\overline{\DD }$. In particular, this shows that all
the polynomials $P(z)$ with $a_1>0$ are contained in a single
open component of $\mathcal{P}^n_{loc}$.


Property (i) easily follows from the continuity arguments and
(\ref{www}).

To prove (ii) we suppose that $P\in \partial \mathcal{P}^n_{loc}$.
Then
$$
\min_{z\in\overline{\DD}}|P'(z)|=0
$$
and it follows from (i) that there is a root $\z_k$ of $P'(z)$
such that $|\z_k|=1$. On the other hand, the coefficients of $P$
are real and it follows that $\overline{\z_k}=1/\z_k$ is a root of
$P'$ as well. But this means that $P'(z)$ and $P'^*(z)$ has a
common root and by the characteristic property of the resultant
the latter is equivalent to $\Res(P',P'^*)=0$.
\end{proof}

\begin{proof}[Proof of Theorem~\ref{ttt}]
Let us consider a real-valued continuous function
$$
f(a)=\Res(S(a),S(a)^*):\mathcal{P}^n_{loc}\to \R{}
$$
where the star is used for the reciprocal polynomial.

We note that $f$ does not change sign on $\mathcal{P}^n_{loc}$.
Indeed, given an arbitrary $P(z)\in \mathcal{P}^n_{loc}$ we have
that all the roots $\z_k$ of $P'(z)$ are outside  $\overline{\DD
}$. Thus,
$$
\qquad |\z_i\z_j|>1, \quad \forall\; i,j\leq n-1.
$$
and by (\ref{A3}) $f(a)\ne 0$. The last inequality together with
(\ref{sign}) implies the claimed property. It easily follows from
the normalization  $a_1>0$ and (\ref{rez_defA}) that $f>0$ on
$\mathcal{P}^n_{loc}$.

Hence, $\mathcal{P}^n_{loc}\subset \Lambda$ for certain open
component $\Lambda$ of $f> 0$. On the other hand, by property (ii)
in Lemma~\ref{lem:com} we have $f(a)=0$ for all $a\in\partial
\mathcal{P}^n_{loc}$. Then it follows by (\ref{sign}) that
$\Lambda=\mathcal{P}^{n}_{loc}$.

Thus, we conclude that $\mathcal{P}^n_{loc}$ coincides with a
certain open component of
$$
\R{n}\setminus\ker W_{n-1}=\R{n}\setminus S^{-1}(\ker f).
$$

To finish the proof we have only to check that the three algebraic
surfaces in the statement of Theorem~\ref{ttt} have nonempty
intersection with the boundary components of $\mathcal{P}^n_{loc}$
(for $n\geq 3$). Indeed, we notice that the hyperplanes
$\Pi^{\pm}$ in (\ref{Pi}) correspond to the polynomials
$P\in\partial \mathcal{P}^n_{loc}$ which have their critical
points on the real axis: $P'(\pm1)=0$. On the other hand,
$\mathcal{A}$ in (\ref{Pi1}) represents the component of $\partial
\mathcal{P}^n_{loc}$ which consists of the polynomials with the
complex roots $\z\not\in\R{}$, $|\z|=1$, $P'(\z)=0$.
\end{proof}

\end{document}